\documentclass{article}
\usepackage{amsmath,amsfonts,amssymb}

\DeclareMathSymbol\nullset{\mathord}{AMSb}{"3F}

\begin{document}

\title{The central polynomials of the finite dimensional unitary and nonunitary Grassmann algebras}

\author{C. Bekh-Ochir and S. A. Rankin}
\maketitle

\begin{abstract}
We describe the $T$-space of central polynomials for both the unitary 
and the nonunitary finite dimensional Grassmann algebra over a field 
of characteristic $p\ne2$ (infinite field in the case of the unitary algebra).
\end{abstract}

%\keywords{central polynomials, $T$-space, $T$-ideal, finite dimensional Grassmann algebra, unitary, nonunitary.}

%\ccode{2000 Mathematics Subject Classification: 22E46, 53C35, 57S20}

\newtheorem{theorem}{Theorem}[section]
\newtheorem{corollary}{Corollary}[section]
\newtheorem{lemma}{Lemma}[section]
\newtheorem{proposition}{Proposition}[section]
{\newtheorem{definition}{Definition}[section]}
\def\proof{\ifdim\lastskip<\smallskipamount\relax\removelastskip
  \vskip\smallskipamount\fi\leavevmode\noindent\hbox to 0pt{\hfil}{\it Proof.}}
\def\strutdepth{\dp\strutbox}
\def\epmarker{\vbox to \strutdepth{\baselineskip\strutdepth\vss\hfill{%
\hbox to 0pt{\hss\vrule height 4pt width 4pt depth 0pt}\null}}}
\def\edproofmarker{\strut\vadjust{\kern-2\strutdepth\epmarker}}
\def\endproof{\edproofmarker\vskip10pt}

%%%%%%%%%%%%%%%%%%%% references from \cite{Ra} %%%%%%%%%%%%%%%%%%
\def\corAltDescriptionOfSGen{2.1} %\label{corollary: alt description of S gen}
\def\corNeedForRGen{2.2} %\label{corollary: need for R gen}
\def\lemForUnitary{3.2} %\label{lemma: for unitary}
\def\thmMainTheoremForNonunitary{2.1} %\label{theorem: main theorem for nonunitary}
\def\corRepOfK1x{3.1}   %\label{corollary: rep of k1x}
\def\lemMprimeFact{3.3} %\label{lemma: m' fact} 
\def\thmCentralInUnitary{3.1} %\label{theorem: central in unitary} 
%%%%%%%%%%%%%%%%%%%%%%%555
\newcounter{parts}
\def\set#1\endset{\{\,#1\,\}}
\def\gen#1{\left<{#1}\right>}
%\newsymbol\nullset 203F  % replaces \emptyset
\def\rest#1{{}_{{}_{#1}}}
\def\com#1,#2{[{#1},{#2}]}
\def\choice#1,#2{\binom{#1}{#2}}
\def\kx{k\langle X\rangle}
\def\kzerox{k_0\langle X\rangle}
\def\konex{k_1\langle X\rangle}
\def\unitgrass{G}
\def\nonunitgrass{G_0}
\def\finiteunitgrass#1{G(#1)}
\def\finitenonunitgrass#1{G_0(#1)}
\def\siderovset{SS}
\def\finitesiderovsetone#1{SS'(#1)}
\def\finitesiderovsettwo#1{SS''(#1)}
\def\boundedsiderovset{BSS}
\def\bss#1{BSS(#1)}
\def\lend#1{lend(#1)}
\def\lbeg#1{lbeg(#1)}
\def\finitesidset#1{SS(#1)}
\let\cong=\equiv
\def\mod#1{\,\,(\text{mod}\,#1)}

%%%%%%%%%%%%%%%%%%%%%%%%%%%%%%%%% Section 1 %%%%%%%%%%%%%%%%%%%%%%%%%%%%%%%

\section{Introduction and preliminaries}
 Let $k$ be a field and $X$ a countable set, say $X=\set x_i\mid i\ge 1\endset$. Then $\kzerox$ denotes
 the free (nonunitary) associative $k$-algebra over $X$, while $\konex$ denotes the free unitary associative
 $k$-algebra over $X$. 

 Let $H$ denote any associative $k$-algebra. 
 For any $X\subseteq H$, $\gen{X}$ shall denote the linear subspace of $H$ spanned by $X$.
 Any linear subspace of $H$ that is invariant under every endomorphism of $H$ is
 called a $T$-space of $H$, and if a $T$-space happens to also be an ideal of $H$, then it
 is called a $T$-ideal of $H$. For $X\subseteq H$, the smallest $T$-space containing $X$
 shall be denoted by $X^S$, while the smallest $T$-ideal of $H$ that contains $X$ shall be
 denoted by $X^T$. In this article, we shall deal only with $T$-spaces and
 $T$-ideals of $\kzerox$ and $\konex$.

 An element $f\in \kzerox$ is called an {\em identity} of $H$
 if $f$ is in the kernel of every homomorphism from $\kzerox$ to $H$ (from $\konex$ if $H$ is unitary). 
 The set of all identities of $H$ is a $T$-ideal of $\kzerox$ (and of $\konex$ if $H$ is unitary), 
 denoted by $T(H)$. An element $f\in \kzerox$ is called a {\it central polynomial} of $H$ 
 if $f\notin T(H)$ and the image of $f$ under any homomorphism from $\kzerox$ ($\konex$ if $H$ is unitary)
 belongs to $C_{H}$, the centre of $H$.

 Let $\unitgrass$ denote the (countably) infinite dimensional unitary
 Grassmann algebra over $k$, so there exist $e_i\in \nonunitgrass$, $i\ge 1$, 
 such that  for all $i$ and $j$, $e_ie_j=-e_je_i$, $e_i^2=0$, and
 $\mathcal{B}=\set e_{i_1}e_{i_2}\cdots e_{i_n}\mid n\ge 1,
 i_1<i_2<\cdots i_n\endset$, together with $1$, forms a linear basis for
 $G$. The subalgebra of $\unitgrass$ with linear basis $\mathcal{B}$ is
 the infinite dimensional nonunitary Grassmann algebra over $k$, and is
 denoted by $\nonunitgrass$. Then for any positive integer $m$, the unitary
 subalgebra of $\unitgrass$ that is generated by $\set e_1,e_2,\ldots,e_m\endset$,
 is denoted by $\finiteunitgrass{m}$, while the nonunitary subalgebra of
 $\nonunitgrass$ that is generated by the same set is denoted by $\finitenonunitgrass{m}$.

 Evidently, $T(\finiteunitgrass{m})\subseteq T(\finitenonunitgrass{m})$. It is well known
 that $T^{(3)}$, the $T$-ideal of $\konex$ that is generated  by $\com
 {\com x_1,{x_2}},{x_3}$, is contained in $T(\finiteunitgrass{m})$. For convenience, 
 we shall write $\com x_1,{x_2,x_3}$ for  $\com {\com x_1,{x_2}},{x_3}$.

 We shall let $CP(\finiteunitgrass{m})$ and $CP(\finitenonunitgrass{m})$ denote
 the $T$-spaces of $\konex$ and $\kzerox$ that are generated by the 
 central polynomials of $\finiteunitgrass{m}$ and $\finitenonunitgrass{m}$, respectively. 
 Evidently, $T(\finiteunitgrass{m})\subseteq CP(\finiteunitgrass{m})$, $T(\finitenonunitgrass{m})
 \subseteq CP(\finitenonunitgrass{m})$, and $CP(\finiteunitgrass{m})\cap \kzerox\subseteq CP(\finitenonunitgrass{m})$.

 For a field of characteristic 2 or if $m=1$, the unitary finite dimensional Grassmann
 algebra of dimension $m$ (and hence the nonunitary Grassmann algebra of dimension $m$) is commutative,
 and thus the $T$-space of central polynomials for each is $\kzerox$.
 It is for this reason that we restrict our attention to those fields of characteristic $p\ne 2$ and
 only consider Grassmann algebras of finite dimension $m\ge 2$.
 
 In this paper, we use techniques developed in \cite{Ra} to describe 
 the $T$-space of central polynomials for
 both the unitary and the nonunitary finite dimensional Grassmann
 algebra over a field of characteristic $p\ne2$ (except for the case of
 the finite dimensional unitary Grassmann algebra over a finite field).
 As in the infinite dimensional case, the description of the central polynomials relies heavily on
 having complete knowledge of the $T$-ideal of identities of the
 relevant Grassmann algebra, due to A. Giambruno and P.
 Koshlukov \cite{Gi} for the case of the unitary finite dimensional
 Grassmann algebra over an infinite field of characteristic $p>2$, and 
 to A. Stojanova-Venkova
 \cite{At} for the case of the nonunitary finite dimensional Grassman
 algebra over an arbitrary field. It is of interest to note that
 the same general approach has been used for both the unitary and
 nonunitary, infinite and finite dimensional Grassman algebras.
 
 We complete this section with a brief description of results from the literature
 that will be required in the sequel.

\begin{lemma}[\cite{Ra}, Lemma 1.1]\label{lemma: handy}\  % line 138
 \vskip-1.5\baselineskip\null
 \begin{list}{(\roman{parts})}{\usecounter{parts}}
  \item % (i)
    $\com u,{vw}=\com u,vw+v\com u,w$ for all $u,v,w\in \kzerox$.
  \item % (ii)
    $\com u,{vw}=\com u,vw+\com u,wv+\com v,{\com u,w}$ for all $u,v,w\in \kzerox$.
  \item % (iii)
     $\com u,{\prod_{i=1}^n v_i}\cong \sum_{i=1}^n 
          \bigl(\mkern3mu \com u,{v_i}\prod_{\substack{j=1\\j\ne i}}^n v_j\mkern3mu\bigr)\mod{T^{(3)}}$ for any positive 
	  integer $n$, and any $u,v_1,v_2,\ldots,v_n\in \kzerox$. 
   \item  % (iv)
     $\com u,v\com w,x\cong -\com u,w\com v,x\mod{T^{(3)}}$ for all $u,v,w\in \kzerox$.
   \item  % (v)
    $\com u,v\com u,w\cong 0\mod{T^{(3)}}$ for all $u,v,w\in \kzerox$.
   \item  % (vi)
    $\com u,v uw\cong \com u,vwu\mod{T^{(3)}}$ for all $u,v,w\in \kzerox$.    
  \item % (vii)
   For any $n\ge2$, $x_1^nx_2^n\cong(x_1x_2)^n+\choice n,2\com x_1,{x_2}x_1^{n-1}x_2^{n-1}\mod{T^{(3)}}$.    
 \end{list}
\end{lemma}

The following lemma summarizes discussion found in Siderov \cite{Si}, where $C$ is the linear subspace of $\nonunitgrass$
spanned by $\set \prod_{i=1}^{2n} e_i\mid n\ge 1\endset$, and $H$ is the linear subspace of $\nonunitgrass$ 
spanned by $\set \prod_{i=1}^{2n-1} e_i\mid n\ge 1\endset$.
 
\begin{lemma}\label{lemma: useful}\  % line 425
 \vskip-1.5\baselineskip\null
  \begin{list}{(\roman{parts})}{\usecounter{parts}}
  \item % (i)
    $C=C_{\nonunitgrass}$. % k+ C=C_G in unitary case
  \item % (ii)
    For $h,u\in H$, $hu=-uh$. In particular, $h^2=0$ (since $p\ne2$).
  \item % (iii)
    Let $g\in \nonunitgrass$, so there exist (unique) $c\in C$ and $h\in H$ such that $g=c+h$. For any %% g=\alpha+c+h \alpha\in k in unitary case
  positive integer $n$, $g^n=c^n+nc^{n-1}h$. 
  \item % (iv)
   For $g\in \nonunitgrass$, $g^p=0$.
  \item % (v)
   Let $c_1,c_2\in C$ and $h_1,h_2\in H$, and set $g_1=c_1+h_1$, $g_2=c_2+h_2$. Then
   for any nonnegative integers $m_1,m_2$, $\com g_1,{g_2}g_1^{m_1}g_2^{m_2}=2c_1^{m_1}c_2^{m_2}h_1h_2$
   (where $g_i^0$ and $c_i^0$ are understood to mean that the factors $g_i^0$ and $c_i^0$ are omitted).
  \item % (vi)
   Let $u\in \nonunitgrass$. Then $u^{n+1}=0$, where $n$ is the number of distinct basic product terms in the expression
   for $u$ as a linear combination of elements of $\mathcal{B}$.
 \end{list}
\end{lemma}

\begin{definition}\label{definition: siderov's elements}
 Let $\siderovset$ denote the set of all elements of the form 
 \begin{list}{(\roman{parts})}{\usecounter{parts}}
  \item $\prod_{r=1}^t x_{i_r}^{\alpha_r}$, or
  \item $\prod_{r=1}^s \com x_{j_{2r-1}},{x_{2r}}x_{j_{2r-1}}^{\beta_{2r-1}}x_{j_{2r}}^{\beta_{2r}}$, or
  \item $\bigl(\prod_{r=1}^t x_{i_r}^{\alpha_r}\bigr)\prod_{r=1}^s \com x_{j_{2r-1}},{x_{2r}}x_{j_{2r-1}}^{\beta_{2r-1}}x_{j_{2r}}^{\beta_{2r}}$, 
 \end{list}  
 \noindent where $j_1<j_2<\cdots j_{2s}$, $\beta_i\ge0$ for all $i$, $i_1<i_2<\cdots < i_t$, 
 $\set i_1,\ldots,i_r\endset\cap\set j_1,\ldots,j_{2s}\endset=\nullset$, and $\alpha_i\ge 1$ for all $i$.

 Let $u\in \siderovset$. If $u$ is of the form (i), then the beginning of $u$ is
 $\prod_{r=1}^t x_{i_r}^{\alpha_r}$, the end of $u$ is empty, the length of the beginning of $u$, $\lbeg{u}$, is
 equal to $t$ and the length of the end of $u$, $\lend{u}$, is 0. If $u$ is of the form
 (ii), then we say that the beginning of $u$ is empty, the end of $u$ is $\prod_{r=1}^s \com x_{j_{2r-1}},{x_{2r}}x_{j_{2r-1}}^{\beta_{2r-1}}x_{j_{2r}}^{\beta_{2r}}$,
 and $\lbeg{u}=0$ and $\lend{u}=s$.
 If $u$ is of the form (iii), 
 then we say that the beginning of $u$ is $\prod_{r=1}^t x_{i_r}^{\alpha_r}$, the end of $u$ is
 $\prod_{r=1}^s \com x_{j_{2r-1}},{x_{2r}}x_{j_{2r-1}}^{\beta_{2r-1}}x_{j_{2r}}^{\beta_{2r}}$, and $\lbeg{u}=t$ and
 $\lend{u}=s$.
\end{definition}

In \cite{At}, Venkova introduced a total order on the set $\siderovset$ which was useful in her work on the
identities of the finite dimensional nonunitary Grassmann algebra. 

%% Va's total order
\begin{definition}[Venkova's ordering]\label{definition: total order}
 For $u,v\in \siderovset$, we say that $u>v$ if one of the following requirements holds.
  \begin{list}{(\roman{parts})}{\usecounter{parts}}
  \item $\deg u < \deg v$.
  \item $\deg u = \deg v$ but $\lend{u}<\lend{v}$.
  \item $\deg u = \deg v$ and $\lend{u}=\lend{v}$, but there exists $i \ge 1$ such that 
        $\deg_{x_i} u < \deg_{x_i} v$ and for each $j<i$, $\deg_{x_j} u = \deg_{x_j} v$.
  \item $\deg u = \deg v$, $\lend{u}=\lend{v}$ and for each $i \ge 1$,
        $\deg_{x_i} u = \deg_{x_i} v$, and there exists $j \ge 1$ such that 
        $x_j$ appears in the end of $u$ and in the beginning of $v$, and 
        for each $k<j$, $x_k$ appears in the beginning of $u$ if and only if $x_k$ appears in 
        the beginning of $v$.
  \end{list}
 \end{definition}
 
 It will be helpful to note that if $u>v$ by virtue of condition (iv), then there exists $k>j$ such that $x_k$ is in
 the beginning of $u$ and in the end of $v$.

%%%%%%%%%%%%%%%%%%%%%%%%%%%%%%%%%%%%%%%%%% Section 2 %%%%%%%%%%%%%%%%%%%%%%%%%%%%%%%%%%%%%%%%%%%%%%
\section{The central polynomials of the finite dimensional unitary Grassmann 
                     algebra over an infinite field of characteristic $p\ne2$}  
 In this section, $k$ denotes an infinite field of characteristic $p\ne2$,
 and for any integer $m\ge 2$, $\finiteunitgrass{m}$ denotes the 
 subalgebra of $\unitgrass$ that is generated (as an algebra) by
 $\set e_1,e_2,\ldots,e_m\endset$.

\begin{lemma}\label{lemma: centre of gm}
 Let $m\ge2$ be an integer. Then $C_{\finiteunitgrass{m}}$ is equal 
 to $C_{\unitgrass}\cap \finiteunitgrass{m}$ if $m$ is even, and to 
 $(C_{\unitgrass}\cap \finiteunitgrass{m})+\gen{\prod_{r=1}^m e_r}$ if $m$ is odd.
\end{lemma}

\begin{proof}
 Since $C_{\unitgrass}\cap \finiteunitgrass{m}\subseteq C_{\finiteunitgrass{m}}$,
 and $e_1e_2\cdots e_mg=ge_1\cdots e_m$ for all $g\in \finiteunitgrass{m}$
 (both products are $0$ for $g\in\gen{\mathcal{B}}$), we have $(C_{\unitgrass}\cap 
 \finiteunitgrass{m})+\gen{\prod_{r=1}^m e_r}\subseteq C_{\finiteunitgrass{m}}$.
 Let $g\in C_{\finiteunitgrass{m}}$. Then there exist $\alpha\in k$, $c\in C\cap \finiteunitgrass{m}$, and
 $h\in H\cap \finiteunitgrass{m}$ such that $g=\alpha+c+h$. For any $e_i$ with $1\le i\le m$, we have 
 $ge_i=e_ig$, but $ge_i=\alpha e_i+ce_i+he_i$, while $e_ig=\alpha e_i+ce_i+e_ih$,
 so $e_ih=he_i$ for each $i$ with $1\le i\le m$. However, since $h\in H$, we have
 $e_ih=-he_i$, and thus $2e_ih=0$ for all $i=1,2,\ldots,m$. Since $p\ne2$, we have
 $e_ih=0$ for all $i=1,2,\ldots,m$. Since the elements of $\mathcal{B}\cap \finiteunitgrass{m}$
 are linearly independent in $\finiteunitgrass{m}$, this implies that $h$ is either 0
 or else is a scalar multiple of $e_1e_2\ldots e_m$ (in which case, $m$ is odd),
 which proves the result.
\end{proof}
 
\begin{definition}\label{definition: product of commutators}
 For any positive integer $j$, let $h_j=\prod_{r=1}^j \com x_{2r-1},{x_{2r}}$.
\end{definition}

\begin{lemma}\label{lemma: identities of unitary fin dim grassman}
 Let $m$ be a positive integer, and let $b_m=\lfloor\frac{m}{2}\rfloor+1$. Then 
 $T(\finiteunitgrass{m})$, the $T$-ideal of identities of the unitary Grassmann
 algebra $\finiteunitgrass{m}$, has basis $\set \com x_1,{x_2,x_3}, h_{b_m}\endset$.
\end{lemma}

\proof
 If $p>2$, this is Corollary 8 of \cite{Gi}, so suppose that $p=0$. By
 \cite{D}, Proposition 4.3.3, $T(\finiteunitgrass{m})$ has as basis the
 set of all proper multilinear identities of $\finiteunitgrass{m}$; that is,
 the set of all multilinear products of commutators. Since $T^{(3)}\subseteq T(\finiteunitgrass{m})$,
 it suffices to consider only multilinear products of 2-commutators. Let $r=b_m$ and
 consider $h_r$. We claim that $h_r\in T(\finiteunitgrass{m})$. Evidently, if any $x_i$, 
 $i=1,2,\ldots,2r$, is replaced by an element of $k$, the result is 0. By the
 multilinearity of $h_r$, it suffices to consider $h_r(g_1,g_2,\ldots,g_{2r})$, where
 $g_i=c_i+h_i$, $c_i\in C$ and $h_i\in H$ for each $i=1,2,\ldots,2r$. Then by
 Lemma \ref{lemma: useful} (v), $h_r(g_1,g_2,\ldots,g_{2r})=2^r\prod_{i=1}^{2r}h_i$.
 Since $2r\ge m+1$, it follows that $h_r(g_1,g_2,\ldots,g_{2r})=0$. Since $h_r\in\kzerox$, this establishes that
 $h_r\in T(\finiteunitgrass{m})$. On the other hand, for any $n<r$ (so $2n\le m$), 
 $h_n(e_1,e_2,\ldots,e_{2n})=2^ne_1e_2\cdots e_{2n}\ne 0$ since $p\ne 2$. Thus
 $T(\finiteunitgrass{m})=\set \com x_1,{x_2,x_3},h_r\endset^T$.
\endproof

\begin{definition}\label{definition: support}
 For $u=e_{i_1}e_{i_2}\cdots e_{i_n}\in \mathcal{B}$, let $s(u)=\set e_{i_1},e_{i_2},\ldots,e_{i_n}\endset$, while 
 $s(1)=\nullset$. Then for any
 $g\in \unitgrass$, let $s(g)=\bigcup_{i=1}^m s(g_i)$, where $g=\sum_{i=1}^m a_ig_i$ with $a_i\in k^*$ 
 and $g_i\in \mathcal{B}\cup\set 1\endset$.
 We shall refer to $s(u)$ as the {\em support} of $u$.
\end{definition}

\begin{lemma}\label{lemma: fin unitary m}
 Let 
 \[
  u=\prod_{r=1}^t x_{r}^{\alpha_r}\text{ or }\bigl(\prod_{r=1}^t
 x_{r}^{\alpha_r}\bigr)\prod_{r=1}^s \com x_{t+2r-1},{x_{t+2r}}
 x_{t+2r-1}^{\alpha_{t+2r-1}}x_{t+2r}^{\alpha_{t+2r}},
 \]
 where $t\ge 1$,
 $s\ge 1$, $\alpha_r\ge 1$ for $1\le r\le t$, and if $u$ is of the latter type, 
 $\alpha_r\ge 0$ for $t+1\le r\le t+2s$. Then  $u\notin CP(\finiteunitgrass{m})$ 
 if and only if $2\lend{u}\le m-2$ and, if $p>2$, there exists $a$ with $1\le a\le t$ 
 and $\alpha_a\not\cong0\mod{p}$.
\end{lemma}

\begin{proof}
 Suppose first that $u\notin CP(\finiteunitgrass{m})$. If $p>2$, then by \cite{Ra}, Lemma
 \lemForUnitary\ and the fact that $CP(\unitgrass)\subseteq
 CP(\finiteunitgrass{m})$, there must exist $a$ with $1\le a\le t$ and
 $\alpha_a\not\cong0\mod{p}$. If $u$ is of the first type, then 
 $\lend{u}=0$ and since $m-2\ge0$, there is nothing more to show. Thus it suffices to
 consider only $u$ of the latter type. Since $u\notin CP(\finiteunitgrass{m})$, there
 exist $g_1,g_2,\ldots,g_n\in\finiteunitgrass{m}$, where $n=t+2s$, such that
 $u(g_1,g_2,\ldots,g_n)\notin C_{\finiteunitgrass{m}}$. 
 For each
 $r=1,\ldots,n$, let $g_r=d_r+h_r$, where $d_r=a_r+c_r$, $a_r\in k$,
 $c_r\in C_m$, $h_r\in H_m$. Then $u(g_1,g_2,\ldots,g_n)$ is a linear
 combination of terms of the form $d_{i_1}d_{i_2}\cdots
 d_{i_e}h_{i_{e+1}}\cdots h_{i_{t}}c$ (including the possibility that no
 $d_{i_j}$'s appear, which we think of as $e=0$, or no $h_{i_j}$'s
 appear, which we think of as $e=t$), where $\set i_1,i_2,\ldots,i_t\endset
 =\set 1,2,\ldots,t\endset$, and $c=\prod_{r=1}^s \com
 g_{t+2r-1},{g_{t+2r}}g_{t+2r-1}^{\alpha_{t+2r-1}}g_{t+2r}^{\alpha_{t+2r}}$.
 It follows from Lemma \ref{lemma: useful} (v) that $c\in
 C_{\finiteunitgrass{m}}$. Since $u(g_1,g_2,\ldots,g_n)\notin C_{\finiteunitgrass{m}}$, at
 least one of these terms must not belong to $C_{\finiteunitgrass{m}}$,
 and so there exists $e$ such that the support sets of $c$,
 $d_{i_1},\ldots,d_{i_e},h_{i_{e+1}},\ldots,h_{i_{t}}$ are  pairwise
 disjoint, $t-e$ is odd, and, by Lemma \ref{lemma: centre of gm}, if $m$
 is odd, the union of the supports sets must be properly contained in
 $\set e_1,e_2,\ldots,e_m\endset$. Since $|s(d_{i_r})|\ge0$ and
 $|s(h_{i_r})|\ge1$, we have 
 \[
  |s(c)|+1\le |s(c)|+t-e\le |s(c)|+\sum_{r=1}^e|s(d_{i_r})| 
  + \sum_{r=e+1}^t |s(h_{i_r})|\le \begin{cases} m &\text{$m$ even}\\
                                        m-1 & \text{$m$ odd.}\end{cases}
 \]

 By Lemma \ref{lemma: useful} (v), $c$ is a scalar multiple (necessarily
 nonzero) of a product of terms of the form
 $c_{j_{2r-1}}^{\alpha_{t+2r-1}}c_{j_{2r}}^{\alpha_{t+2r}}h_{t+2r-1}h_{t+2r}$,
 so $|s(c)|\ge 2s$. Thus we find that if $m$ is even, then  $1+2s\le m$,
 or $2s\le m-1$, while if $m$ is odd, then $1+2s\le m-1$, or $2s\le
 m-2$. For $m$ even, $2s\le m-1$ is equivalent to $2s\le m-2$, so these
 two cases can be  consolidated into simply $2s\le m-2$, as required. 
 
 Conversely, we prove that if $2\lend{u}\le m-2$, and, if $p>2$, there exists 
 $a$ with $1\le a\le t$ and $\alpha_a\ne 0\mod{p}$, then $u\notin
 C(\finiteunitgrass{m})$. In the characteristic zero case, we take $a=1$. Now set $g_a=1+e_a$, and $g_i=1$ for all $i$ with
 $1\le i\le t$, $i\ne a$. If $u$ is of the latter type, set 
 $g_{t+r}=1+e_{1+r}$ for $r$ with $1\le r\le 2s$. Then the odd part of
 $u(g_1,g_2,\ldots,g_n)=(1+\alpha_a e_a)
 \prod_{r=1}^{1+s}2e_{1+2r-1}e_{1+2r}$ is equal to 
 $\alpha_a e_a\prod_{r=1}^{1+s}2e_{1+2r-1}e_{1+2r}$, and this is nonzero since $\alpha_a\ne 0$ and
 $p\ne 2$. Since the support
 of $e_a\prod_{r=1}^{1+s}e_{1+2r-1}e_{1+2r}$ has size $2s+1\le m-1$,
 it follows that this element is not in  $C_{\finiteunitgrass{m}}$. Thus
 $u(g_1,g_2,\ldots,g_n)\notin C_{\finiteunitgrass{m}}$, and so $u\notin
 CP(\finiteunitgrass{m})$.
\end{proof} 

\begin{definition}\label{definition: def of Mmnprime}
 For positive integers $t,n,r_1,r_2,\ldots,r_n$ with $t\le n$ and, if $p>2$,
 $r_t\not\cong 0\mod{p}$, let $M'_{t,n}(r_1,r_2,\ldots,r_n)$ denote the
 set of all elements $u\in\siderovset$ with $\lbeg{u}>0$, 
 subject to the requirements that $x_t$ appears in the beginning of $u$
 and $\deg_{x_i}(u) = r_i$ for every $i$.
\end{definition}

\begin{lemma}[\cite{Ra}, Lemma \lemMprimeFact]\label{lemma: m' fact}
  Let $t,n,r_1,r_2,\ldots,r_n$ be positive integers with $t\le n$ and, if $p>2$, $r_t\not\cong 0\mod{p}$.
  Then for any $u \in M'_{t,n}(r_1,\ldots,r_n)$, there exist
  $g_1,g_2,\ldots,g_n$ in $\finiteunitgrass{z}$, where $z=2(\deg(u)-\lend{u})-1$, such that the following hold:
   \begin{list}{(\roman{parts})}{\usecounter{parts}}
   \item $u(g_1,g_2,\ldots,g_n)$ has nonzero odd part in $\gen{\prod_{r=1}^z e_i}$;
   \item for any $v\in M'_{t,n}(r_1,\ldots,r_n)$ with $u>v$, $v(g_1,g_2,\ldots,g_n)=0$;
  \end{list}
\end{lemma}

\begin{corollary}\label{corollary: m result for finite unitary}
  For any positive integers $t,n,r_1,r_2,\ldots,r_n$ with $t\le n$ 
  and, if $p>2$, $r_t\not\cong 0\mod{p}$, 
  \[
   \gen{M'_{t,n}(r_1,\ldots,r_n)}\cap CP(\finiteunitgrass{m})=
     \gen{\set u\in M'_{t,n}(r_1,\ldots,r_n)\mid 2\lend{u}\ge m-1\mkern-3mu\endset}.
  \] 
\end{corollary}

\begin{proof}
 By Lemma \ref{lemma: fin unitary m}, $\gen{\set u\in M'_{t,n}(r_1,\ldots,r_n)\mid 2\lend{u}\ge m-1\endset}$
 is contained in $\gen{M'_{t,n}(r_1,\ldots,r_n)}\cap CP(\finiteunitgrass{m})$.
 Let 
 $$
  u\in \gen{M'_{t,n}(r_1,\ldots,r_n)}\cap CP(\finiteunitgrass{m}),
 $$
 with $u\ne 0$. Then $u$ is a linear combination of some $u_1>u_2>\cdots >u_l\in
 M'_{t,n}(r_1,\ldots,r_n)$. Suppose that $u_1\notin CP(\finiteunitgrass{m})$. Then
 by Lemma \ref{lemma: fin unitary m}, $2\lend{u_1}\le m-2$. But then by Lemma \ref{lemma: m' fact}, there exist $g_1,g_2,\ldots,g_n$ in $\finiteunitgrass{z}$, 
 the subalgebra of $\unitgrass$ that is generated by 
 $\set e_1,e_2,\ldots,e_{z}\endset$, where $z=2\lend{u}+1\le m-1$, such that 
 $u_1(g_1,g_2,\ldots,g_n)$ has nonzero odd part in $\gen{\prod_{r=1}^z e_i}$ and 
 $u(g_1,g_2,\ldots,g_n)=u_1(g_1,\ldots,g_n)$. By Lemma \ref{lemma: centre of gm},
 $u(g_1,g_2,\ldots,g_n)\notin C_{\finiteunitgrass{m}}$, and this contradicts the
 fact that $u\in CP(\finiteunitgrass{m})$. Thus $u_1\in CP(\finiteunitgrass{m})$, 
 and then by Lemma \ref{lemma: fin unitary m}, $2 \lend{u_1}\ge m-1$. Since all elements
 of $M'_{t,n}(r_1,\ldots,r_n)$ have the same degree, it follows that
 $\lend{u_i}\ge \lend{u_1}$ for every $i$, and so $2\lend{u_i}\ge m-1$
 for every $i$, as required.
\end{proof}

\begin{corollary}\label{corollary: main fact for finite unitary}
 For any positive integers $t,n,r_1,r_2,\ldots,r_n$ with $t\le n$ and, if $p>2$,
 $r_t\not\cong 0\mod{p}$,  
 \[
  \gen{M'_{t,n}(r_1,\ldots,r_n)}\cap
 CP(\finiteunitgrass{m})\subseteq \set x_{2b_m-1}h_{b_m-1}\endset^S\subseteq CP(\finiteunitgrass{m}).
 \]
\end{corollary}

\begin{proof}
 It is immediate that each element of 
 \[\gen{\set u\in M'_{t,n}(r_1,\ldots,r_n)\mid 2\lend{u}\ge m-1\endset}
 \]
 is contained
 in $\set x_{2b_m-1}h_{b_m-1}\endset^S$, so the first containment
 follows from Corollary \ref{corollary: m result for finite unitary}. That $\set x_{2b_m-1}h_{b_m-1}\endset^S\subseteq CP(\finiteunitgrass{m})$
 follows from Lemma \ref{lemma: fin unitary m}.
\end{proof}

\begin{definition}\label{definition: def of s1m}
 Let $S_1=\set \com x_1,{x_2}\endset^S$ if $p=0$, otherwise, for each
 $n\ge1$, let $w_n=\prod_{k=1}^n\com {x_{2k-1}},{x_{2k}}x_{2k-1}^{p-1}x_{2k}^{p-1}$, and then set
 \[
   S_1=\set \com x_1,{x_2},x_1^p\endset^S+\set x_{2k+1}^pw_k \mid k\ge1\endset^S,
 \]    
 and then for any positive integer $m$, let $S_1(m)=S_1+\set x_{2b_m-1}h_{b_m-1}\endset^S$.
\end{definition}

\begin{definition}\label{definition: def of R1}
 Let $R_1$ be the subspace of $\konex$ that is spanned by  $1$ and
 \begin{align*}
  \{\, u\in \siderovset \mid \lbeg{u}>0,\ &\text{and, if $p>2$, for some $x_i$ in the}\\
   &\text{beginning of $u$, $\deg_{x_i}(u)\not\cong 0\mod{p}$}\,\}.
  \end{align*}
\end{definition} 

\begin{lemma}\label{lemma: finite unitary k1x}
 For any positive integer $m$, $\konex=R_1+S_1(m)+T(G_m)$. 
\end{lemma}

\begin{proof}
 This is immediate from Lemma \ref{lemma: identities of unitary fin dim grassman}
 and \cite{Ra}, Corollary \corRepOfK1x.
\end{proof}

\begin{theorem}\label{theorem: central polys for finite unitary grassman}
 For any positive integer $m$, $CP(\finiteunitgrass{m})=S_1(m)+T(G_m)$.
\end{theorem}

\begin{proof}
 Let $U_1=S_1(m)+T(G_m)$. Then Corollary \ref{corollary: main fact for
 finite unitary} and \cite{Ra}, Theorem \thmCentralInUnitary, together
 with the fact that $CP(\unitgrass)\subseteq CP({\finiteunitgrass{m}})$,
 imply that $S_1(m)\subseteq CP(\finiteunitgrass{m})$ and thus
 $U_1\subseteq CP(\finiteunitgrass{m})$. Suppose that
 $CP(\finiteunitgrass{m})-U_1\ne\nullset$. Since $k$ is infinite, if
 every multihomogeneous element of $CP(\finiteunitgrass{m})$ belonged to
 $U_1$, then $CP(\finiteunitgrass{m})\subseteq U_1$, so there exists a
 multihomogeneous $f\in CP(\finiteunitgrass{m})-U_1$,  say of type
 $(r_1,r_2,\ldots,r_n)$. By Lemma \ref{lemma: finite unitary k1x},
 $f\cong\sum_{i=1}^b\alpha_i u_i\mod{U_1}$, where $u_i\in R_1$ for each
 $i$.  Just as in the proof of \cite{Ra}, Theorem \thmCentralInUnitary, 
 we may assume that each $u_i$ is multihomogeneous of the same
 type  as $f$.

 Observe that for every $j$, $\lbeg{u_j}>0$, and, if $p>2$, there must exist 
 an index $i$ such that $r_i\not\cong 0\mod{p}$ and
 $x_i$ is in the beginning of $u_j$ for some $j$, because
 otherwise $u_j\in U_1$ for all $j$ and so $f\in U_1$, which is not the case. 
 Let $d$ denote the maximum index such that 
 $x_d$ is in the beginning of $u_j$ for some $j$, and, if $p>2$,
 $r_d\not\cong 0\mod{p}$. Suppose that for some index $j$, 
 $x_d$ appears in the end of $u_j$. Without loss of generality, we may
 assume that  $u_j=x_{i_1}^{r_{i_1}}x_{i_2}^{r_{i_2}}\cdots
 x_{i_t}^{r_{i_t}}\prod_{k=1}^s\bigl(\com x_{j_{2k-1}},{x_{j_{2k}}}
 x_{j_{2k-1}}^{\beta_{2k-1}}x_{j_{2k}}^{\beta_{2k}}\bigr)$,
 where $\beta_i=r_{j_i}-1$ for each $i$,  and $d=j_{2l-1}$ for some $l$
 with $1\le l\le s$.
 By \cite{Ra}, Corollary \corNeedForRGen\ (since $S\subseteq S_1$),
 \begin{align*}
   u_j&\overset{U_1}{\cong} r_d^{-1}\com x_{j_{2l}},{x_{i_1}^{r_{i_1}}x_{i_2}^{r_{i_2}}\cdots x_{i_t}^{r_{i_t}}}x_d^{r_d}x_{j_{2l}}^{\beta_{2l}}
     \prod_{\substack{ k=1\\k\ne l}}^s\bigl(\com x_{j_{2k-1}},{x_{j_{2k}}}x_{j_{2k-1}}^{\beta_{2k-1}}x_{j_{2k}}^{\beta_{2k}}\bigr)\\
     &\overset{U_1}{\cong}\sum_{a=1}^t r_d^{-1}r_{i_a}\com x_{j_{2l}},{x_{i_a}}x_{i_a}^{r_{i_a}-1}x_{j_{2l}}^{\beta_{2l}}\bigl(\prod_{\substack{ h=1\\h\ne a}}^t x_{i_h}^{r_{i_h}}\bigr)x_d^{r_d}  
     \prod_{\substack{ k=1\\k\ne l}}^s\bigl(\com x_{j_{2k-1}},{x_{j_{2k}}}x_{j_{2k-1}}^{\beta_{2k-1}}x_{j_{2k}}^{\beta_{2k}}\bigr).\\
     &\overset{T^{(3)}}{\cong}\sum_{a=1}^t r_d^{-1}r_{i_{a}}\bigl(\prod_{\substack{ h=1\\h\ne a}}^t x_{i_h}^{r_{i_h}}\bigr)x_d^{r_d}  
     \com x_{j_{2l}},{x_{i_a}}x_{i_a}^{r_{i_a}-1}x_{j_{2l}}^{\beta_{2l}}\prod_{\substack{ k=1\\k\ne l}}^s\bigl(\com x_{j_{2k-1}},{x_{j_{2k}}}x_{j_{2k-1}}^{\beta_{2k-1}}x_{j_{2k}}^{\beta_{2k}}\bigr).\\
 \end{align*}
 If $p>2$, then since $\com x_k^p,{x_d}\cong p\com x_k,{x_d}x_k^{p-1}=0\mod{T^{(3)}}$,
 it follows that $x_d^{r_d}$ can be moved, modulo $T^{(3)}$, to the left of
 each $x_{i_a}^{r_{i_a}}$ in the beginning of $u_j$ for which $m<i_a$
 (by choice of $m$, $r_{i_a}\cong 0\mod{p}$ if $m<i_a$). Furthermore, by
 application of Lemma \ref{lemma: handy} (iv) and (vi), the end of each 
 summand may be manipulated modulo $T^{(3)}$ so as to present $u_j$ as a 
 sum of multihomogeneous elements of type $(r_1,r_2,\ldots,r_n)$, each 
 in $R_1-U_1$.  We may therefore assume that $x_d$ appears in the 
 beginning of each $u_j$. But then $\sum_{i=1}^b \alpha_i u_i\in 
 M'_{m,n}(r_1,\ldots,r_n)$. Since $f\cong\sum_{i=1}^b \alpha_i u_i\mod{U_1}$, 
 it follows that $\sum_{i=1}^b \alpha_i u_i\in CP(\finiteunitgrass{m})$, so by
 Corollary \ref{corollary: main fact for finite unitary},
 $\sum_{i=1}^b \alpha_iu_i\in \set x_{2b_m-1}h_{b_m-1}\endset^S\subseteq U_1$, 
 and thus $f\in U_1$. This contradicts our choice of $f$, and
 so we must have $CP(\finiteunitgrass{m})-U_1=\nullset$, as required.
\end{proof}

\begin{corollary}\label{corollary: t space basis for unitary cp}
 Let $m\ge2$. If $k$ is a field of characteristic zero, then 
 \[
   CP(\finiteunitgrass{m})=
 \set \com x_1,{x_2},\com x_1,{x_2}\com x_3,{x_4}, x_{2b_m-1}h_{b_m-1}\endset^S,
 \]
 where $b_m=\left\lfloor\frac{m}{2}\right\rfloor+1$, and $h_{b_{m}}=\prod_{r=1}^{b_m} \com x_{2r-1},{x_{2r}}$,
 while if $k$ is an infinite field of characteristic $p>2$, then 
 \begin{align*}
  CP(\finiteunitgrass{m})=
 \{\,  \com x_1,{x_2},\com x_1,{x_2}&\com x_3,{x_4}, x_1^p, x_{2b_m-1}h_{b_m-1}\,\}^S\\
 &+ \{\, x_{2k+1}^pw_k\mid 1\le k\le b_m-1\,\}^S,
 \end{align*}
 where for any $k\ge1$, $w_k=\prod_{r=1}^k\com {x_{2r-1}},{x_{2r}}x_{2r-1}^{p-1}x_{2r}^{p-1}$.
\end{corollary}

%%%%%%%%%%%%%%%%%%%%%%%%%%%%%%%%%%%%%%%%%%%% Section 3 %%%%%%%%%%%%%%%%%%%%%%%%%%%%%%%%%%%%%%%%%%%%%%%%%%%%%%%%%%%%%%%%5555

\section{The central polynomials of the finite dimensional nonunitary 
                 Grassmann algebra over a field of characteristic $p\ne2$}  
 In this section, $k$ denotes a field of characteristic $p\ne2$, and for
 any positive integer $m$, $\finitenonunitgrass{m}$ denotes the 
 subalgebra of $\nonunitgrass$ that is generated (as an algebra) by
 $\set e_1,e_2,\ldots,e_m\endset$.

The proof of the following lemma is very similar to that of Lemma \ref{lemma: centre of gm}
and will be omitted.

\begin{lemma}\label{lemma: centre of finite nonunitary grassman}
 Let $m$ be a positive integer. Then $C_{\finitenonunitgrass{m}}$ is
 equal to $C_{\nonunitgrass} \cap \finitenonunitgrass{m}$ if $m$ is
 even, and to $(C_{\nonunitgrass} \cap \finitenonunitgrass{m})
 +\gen{\prod_{r=1}^m e_i}$ if $m$ is odd.
\end{lemma} 

%\begin{proof}
% Since $C_{\nonunitgrass}\cap \finitenonunitgrass{m}\subseteq C_{\finitenonunitgrass{m}}$,
% and $e_1e_2\cdots e_mg=0=ge_1\cdots e_m$ for all $g\in \finitenonunitgrass{m}$,
% we have $(C_{\nonunitgrass}\cap 
% \finitenonunitgrass{m})+\gen{\prod_{r=1}^m e_r}\subseteq C_{\finitenonunitgrass{m}}$.
% Let $g\in C_{\finitenonunitgrass{m}}$. Then there exist $c\in C\cap \unitgrass{m}$ and
% $h\in H\cap \nonunitgrass{m}$ such that $g=c+h$. For any $e_i$ with $1\le i\le m$, we have 
% $ge_i=e_ig$, but $ge_i=ce_i+he_i$, while $e_ig=ce_i+e_ih$,
% so $e_ih=he_i$ for each $i$ with $1\le i\le m$. However, since $h\in H$, we have
% $e_ih=-he_i$, and thus $2e_ih=0$ for all $i=1,2,\ldots,m$. Since $p>2$, we have
% $e_ih=0$ for all $i=1,2,\ldots,m$. Since the elements of $\mathcal{B}\cap \finitenonunitgrass{m}$
% are linearly independent in $\finitenonunitgrass{m}$, this implies that $h$ is either 0
% or else is a scalar multiple of $e_1e_2\ldots e_m$ (in which case, $m$ is odd),
% which proves the result.
%\end{proof}

 Let $x_1\circ x_2=x_1x_2+x_2x_1$, and for any $n\ge3$, $x_1\circ
 x_2\circ \cdots\circ x_n=(x_1\circ x_2\circ\cdots\circ x_{n-1})\circ
 x_n$. We remark that $\circ$ is associative modulo $T^{(3)}$.
 
\begin{lemma}[\cite{At}, Theorem 1]\label{lemma: identities of nonunitary fin dim grassman}
 For any positive integer $m$, $T(\finiteunitgrass{m})$, the $T$-ideal of
 identities of the nonunitary Grassmann algebra $\finiteunitgrass{m}$, is 
 generated as a $T$-ideal by:
  \begin{list}{(\roman{parts})}{\usecounter{parts}} 
  \item $x_1^2$ and $x_1x_2\cdots x_{m+1}$ if $p=2$;
  \item $x_1^p$, $\com x_1,{x_2,x_3},x_1\circ x_2\circ \cdots\circ x_{\frac{m}{2}+1}$, 
        if $p>2$ and $m$ is even;
  \item $x_1^p$, $\com x_1,{x_2,x_3},(x_1\circ x_2\circ \cdots 
        \circ x_{\frac{m+1}{2}})x_{\frac{m+3}{2}}, x_{\frac{m+3}{2}}(x_1
        \circ x_2\circ \cdots\circ x_{\frac{m+1}{2}})$, and, if $2p-1$ 
        divides $m+1$, $\prod_{r=1}^{\frac{m+1}{2(2p-1)}}
        \com x_{2r-1},{x_{2r}}x_{2r-1}^{p-1}x_{2r}^{p-1}$,
        if $p>2$ and $m$ is odd;
  \item  $\com x_1,{x_2,x_3}$ and $x_1\circ\cdots\circ x_{\frac{m}{2}+1}$ when $m$ is
 even and $p=0$, and by $\com x_1,{x_2,x_3}$, $(x_1\circ x_2\circ \cdots\circ
 x_{\frac{m+1}{2}})x_{\frac{m+3}{2}}$, and $x_{\frac{m+3}{2}}(x_1\circ
 x_2\circ \cdots \circ x_{\frac{m+1}{2}})$ when $m$ is odd and $p=0$.        
 \end{list}  
\end{lemma}

 The following corollary is a direct consequence of Lemma \ref{lemma:
 identities of nonunitary fin dim grassman} and the fact that if
 $f\in\kzerox$ is such that  for some variable $x$ that does not appear
 in $f$, we have $xf, fx\in T(\finitenonunitgrass{m})$, then $f\in
 CP(\finitenonunitgrass{m})$.

\begin{corollary}\label{corollary: new finite nonunitary central guy}
 For any odd positive integer $m$, $x_1\circ x_2\circ \cdots \circ 
 x_{\frac{m+1}{2}}\in CP(\finitenonunitgrass{m})$.
\end{corollary}

\begin{definition}\label{definition: def of S(m) and U(m) for finite nonunitary}
 If $p=0$, let $S=\set\com x_1,{x_2}\endset^S$, while if $p>2$, then for each $n\ge1$, let 
 $w_n=\prod_{k=1}^n\com {x_{2k-1}},{x_{2k}}x_{2k-1}^{p-1}x_{2k}^{p-1}$ and set
 $S=\set \com x_1,{x_2}\endset^S+\set w_n\mid \ n\ge 1\endset^S$.
 Then for each positive integer $m$, let
  \[
   S(m)=
     \begin{cases} S+\set x_1\circ\cdots \circ x_{\frac{m+1}{2}}\endset^S&
              \quad\text{if $m$ is odd}\\
            S&\quad\text{if $m$ is even,}
     \end{cases}
  \]             
  and let $U(m)=S(m)+T(\finitenonunitgrass{m})$.
\end{definition}

 Our objective is to prove that for each positive integer $m$, 
 $CP(\finitenonunitgrass{m})=U(m)$. It will be convenient 
 to have the following notation.

\begin{definition}\label{definition: finite siderov set}
 If $p=0$, let $\boundedsiderovset=\siderovset$, while if $p>2$, let 
 $$
  \boundedsiderovset=\set u\in \siderovset\mid \text{ for any $i$, 
 if $x_i^\alpha$ is a factor of $u$, then $\alpha\le p-1$}\endset.
 $$
 Then for any positive integer $m$, let 
%% corrected definition of \bss{m}
 \begin{align*}
  \bss{m}&=\{\, u\in \boundedsiderovset\mid \deg(u)\le m,\ \deg(u)-\lend{u}\le (m+1)/2,\\ 
  &\hskip60pt\text{ for any $i$, if $x_i^\alpha$ is a factor of $u$, then $\alpha\le \frac{m+1}{2}$}\,\}\\
 \finitesiderovsetone{m}&=\set u\in\bss{m}\mid \lend{u}=0\endset\\
 \finitesiderovsettwo{m}&=\set u\in \bss{m}\mid \lbeg{u}\lend{u}>0 \endset.
% \begin{align*}
% \bss{m}&=\set u\in \boundedsiderovset\mid \deg(u)\le m,\ \deg(u)-\lend{u}\le (m+1)/2\endset,\\
 \end{align*}
 Finally, let $\finitesidset{m}=\finitesiderovsetone{m}\cup
 \finitesiderovsettwo{m}$. If $u\in \bss{m}$ satisfies 
 $\deg(u)-\lend{u}=(m+1)/2$, then $u$ is said to be {\em extremal}, 
 otherwise $u$ is said to be non-extremal. 
\end{definition} 
 
 Note that if $\bss{m}$ contains an extremal element, necessarily $m$ is odd.

\begin{definition}
 A polynomial $f\in \konex$ is said to be {\it  essential in its variables} if every variable 
 that appears in any monomial of $f$ appears in every monomial of $f$.
\end{definition}

\begin{lemma}[\cite{Ra}, Lemma 2.7]\label{lemma: basic essential} % line 114
 Let $V$ be a $T$-space in $\kzerox$.
 \begin{list}{(\roman{parts})}{\usecounter{parts}}
  \item Let $f\in V$. If $f$ is not essential in its variables, then there exist $f_0,f_1\in V$ such
        that $f=f_0+f_1$ and $M(f_0)<M(f)$, $M(f_1)<M(f)$.
  \item Let $E_V=\set f\in V\mid f\text{ is essential in its variables}\endset$. Then $V=\gen{E_V}$; that 
        is, $V$ is the linear span of $E_v$.        
  \end{list}
\end{lemma}

In \cite{At}, Venkova states that she considers only fields of characteristic $p>2$, but
many of her results (with proofs) are valid unchanged for characteristic zero, and others are valid with
obvious modifications for the case of characteristic zero. In fact, at the end of her paper, Venkova states that
the proof of the characteristic zero results is analagous but simpler to that of the characteristic $p>2$ case.
In particular, with the definition of $\boundedsiderovset$
above in the case of characteristic zero, the following two lemmas from \cite{At} are valid in characteristic
zero as well as characteristic $p>2$.

\begin{lemma}[Venkova \cite{At}, Lemma 6]\label{lemma: at, lemma 6}
 Let $m$ be even, and $f$ be essential in $\kzerox$, depending on 
 $x_1,x_2,\ldots,x_n$. Then there exist
 $u_r\in \bss{m}$ and $a_r\in k$ for $r\in F$, $F$ a finite set of positive 
 integers, such that 
 $$
  f\cong \sum_{r\in F} \alpha_r u_r\mod{T(\finitenonunitgrass{m})}.
 $$
%% unnecessary in view of corrected definition of \bss{m}
% and for each $r\in F$, the following hold.
%  \begin{list}{(\roman{parts})}{\usecounter{parts}}
%  \item % (i)
%   If $\lbeg{u_r}>0$, then $1\le\alpha_{i_l}<\frac{m}{2}+1$ for each 
%      $l$ with $1\le l\le \lbeg{u_r}$.
%  \item % (ii)
%   If $\lend{u_r}>0$, then $0\le\beta_{j_l}<\frac{m}{2}+1$ for each 
%      $l$ with $1\le l\le 2\lend{u_r}$.
% \end{list}
\end{lemma}   

\begin{lemma}[Venkova, \cite{At}, Lemma 11]\label{lemma: at, lemma 11}
 Let $m$ be odd, and $f$ be essential in $\kzerox$, depending on
 $x_1,x_2,\ldots,x_n$. Then there exist $u_r\in \bss{m}$ and $a_r\in k$
 for $r\in F$, $F$ a finite set of positive integers, such that 
 $$
  f\cong
 \sum_{r\in F} \alpha_r u_r\mod{T(\finitenonunitgrass{m})},
 $$
 where for each $r\in F$, 
% now covered by the revised definition of \bss{m}
%  \begin{list}{(\roman{parts})}{\usecounter{parts}}
%  \item % (i)
%   If $\lbeg{u_r}>0$, then $1\le\alpha_{i_l}<\frac{m+3}{2}$ for each 
%           $l$ with $1\le l\le \lbeg{u_r}$.
%  \item % (ii)
%   If $\lend{u_r}>0$, then $0\le\beta_{j_l}<\frac{m+3}{2}$ for each 
%           $l$ with $1\le l\le 2\lend{u_r}$.
%  \item % (iii)
   if $u_r$ is extremal, then $\lbeg{u_r}>0$ and $i_1$, the index of the
   first variable in the beginning of $u_r$, is less than $j_r$ for any variable
   $x_{j_r}$ that appears in the end of $u_r$ with degree $p$.
% \end{list}
\end{lemma}   

\setbox0=\hbox{$\scriptstyle\frac{n}{2}$}
\dp0=2.5pt
\ht0=0pt
\wd0=5pt
\begin{lemma}\label{lemma: deal with Venkova's guy}
 Let $n$ be a positive integer, and let $J_n=\set 1,2,\ldots,n\endset$. Then
 \[
  x_1\circ x_2\circ\cdots\circ x_n\cong\sum_{s=0}^{\left\lfloor\mkern-5mu\copy0\right\rfloor}
  (-1)^s2^{n-1-s}\mkern-10mu\sum_{\substack{ J\subseteq J_n\\|J|=2s}}\mkern -8mu P_n(J)Q(J)\ \mod{T^{(3)}},
 \]
 where $P_n(J)=1$ if $2s=n$, otherwise $P_n(J)=\prod_{r=1}^{n-2s}
 x_{i_r}$ with $J_n-J=\set i_1,i_2,\ldots,i_{n-2s}\endset$ and
 $i_1<i_2<\cdots <i_{n-2s}$, and $Q(J)=1$ if $J=\nullset$ (that is;
 $s=0$), else $Q(J)=\prod_{r=1}^s\com x_{j_{2r-1}},{x_{j_{2r}}}$, with
 $J=\set j_{1},
 \ldots, j_{2s}\endset$, and $j_1<j_2<\cdots <j_{2s}$.
\end{lemma}

\begin{proof}
 The proof is by induction on $n$, with base case $n=1$. When $n=1$, the
 value of the sum is $P_1(\nullset)=x_1$, as required. Suppose now that
 $n\ge1$ is such that the assertion holds. Then modulo $T^{(3)}$, we have
 \begin{align*}\setbox0=\hbox{$\scriptstyle\frac{n}{2}$}\dp0=2.5pt\ht0=0pt\wd0=5pt
  x_1\circ x_2\circ\cdots \circ x_n\circ x_{n+1}&=2(x_1\circ x_2\circ\cdots \circ x_n)x_{n+1}-\com x_1\circ x_2\circ\cdots \circ x_n,{x_{n+1}}\\
  &\hskip-100pt\cong \setbox0=\hbox{$\scriptstyle\frac{n}{2}$}\dp0=2.5pt\ht0=0pt\wd0=5pt
  2\mkern-3mu\sum_{s=0}^{\left\lfloor\mkern-5mu\copy0\right\rfloor}(-1)^s2^{n-1-s}\mkern-12mu\sum_{\substack{ J\subseteq J_n\\|J|=2s}}
 \mkern -13mu P_n(\mkern-2mu J)x_{n+1}\mkern-2mu Q(\mkern-2mu J)
 -
 \mkern-5mu\sum_{s=0}^{\left\lfloor\mkern-5mu\copy0\right\rfloor}(-1)^s2^{n-1-s}\mkern-13mu\sum_{\substack{ J\subseteq J_n\\|J|=2s}}
 \mkern -8mu \com P_n(\mkern-2mu J),{x_{n+1}}Q(\mkern-2mu J)\\
 &\hskip-100pt\cong \setbox0=\hbox{$\scriptstyle\frac{n}{2}$}\dp0=2.5pt\ht0=0pt\wd0=5pt
 \mkern-2mu\sum_{s=0}^{\left\lfloor\mkern-5mu\copy0\right\rfloor}(-1)^s2^{(n+1)-1-s}\mkern-15mu\sum_{\substack{ J\subseteq J_{n+1}\\|J|=2s\\x_{n+1}\notin J}}
 \mkern -12mu P_{n+1}(\mkern-2mu J)Q(\mkern-2mu J)
 -\mkern-2mu
 \sum_{s=0}^{\left\lfloor\mkern-5mu\copy0\right\rfloor}(-1)^s2^{n-1-s}\mkern-25mu\sum_{\substack{ J\subseteq J_{n+1}\\|J|=2(s+1)\\x_{n+1}\in J}}
 \mkern -18mu P_{n+1}(J)Q(J).
 \end{align*}

 Now if $n$ is even, then
 $\left\lfloor\frac{n}{2}\right\rfloor=\left\lfloor\frac{n+1}{2}\right\rfloor$,
 while if $n$ is odd, then
 $\left\lfloor\frac{n+1}{2}\right\rfloor=\left\lfloor\frac{n}{2}\right\rfloor+1$
 and it is not possible to have $J\subseteq J_{n+1}$  with $|J|=n+1$ and
 $x_{n+1}\notin J$. Thus in the first sum above, we may change the upper
 summation variable limit from  $\left\lfloor\frac{n}{2}\right\rfloor$
 to $\left\lfloor\frac{n+1}{2}\right\rfloor$ without changing the value
 of the sum. Furthermore, 
 \begin{align*}
   \setbox0=\hbox{$\scriptstyle\frac{n}{2}$}\dp0=2.5pt\ht0=0pt\wd0=5pt
   -\sum_{s=0}^{\left\lfloor\mkern-5mu\copy0\right\rfloor}(-1)^s2^{n-1-s}\mkern-20mu\sum_{\substack{ J\subseteq J_{n+1}\\|J|=2(s+1)\\x_{n+1}\in J}}
   \mkern -20mu P_{n+1}(J)Q(J)&=
   \mkern-3mu\setbox0=\hbox{$\scriptstyle\frac{n}{2}$}\dp0=2.5pt\ht0=0pt\wd0=5pt
   \sum_{s=1}^{\left\lfloor\mkern-5mu\copy0\right\rfloor+1}(-1)^{s}2^{(n+1)-1-s}\mkern-15mu\sum_{\substack{ J\subseteq J_{n+1}\\|J|=2s\\x_{n+1}\in J}}
   \mkern -8mu P_{n+1}(J)Q(J)\\
    &=\setbox0=\hbox{$\scriptstyle\frac{n+1}{2}$}\dp0=2.5pt\ht0=0pt\wd0=13pt%
   \sum_{s=0}^{\left\lfloor\mkern-5mu\copy0\right\rfloor}(-1)^{s}
   2^{(n+1)-1-s}\mkern-15mu\sum_{\substack{ J\subseteq J_{n+1}\\|J|=2s\\x_{n+1}\in J}}
   \mkern -8mu P_{n+1}(J)Q(J)
  \end{align*}
 as there is no choice of $J\subseteq J_{n+1}$ with $|J|=0$ and
 $x_{n+1}\in J$, and if $n$ is even, there is no $J\subseteq J_{n+1}$ of 
 size $n+2$, while if $n$ is odd, then
 $\left\lfloor\frac{n}{2}\right\rfloor+1=\left\lfloor\frac{n+1}{2}\right\rfloor$.   
 It now follows that modulo $T^{(3)}$,  $x_1\circ x_2\circ\cdots \circ
 x_n\circ x_{n+1}$ is congruent to
 \setbox0=\hbox{$\scriptstyle\frac{n+1}{2}$}\dp0=2.5pt\ht0=0pt\wd0=13pt%
 \[
   \sum_{s=0}^{\left\lfloor\copy0\right\rfloor}\mkern-4mu (-1)^s2^{(n+1)-1-s}\mkern-15mu\sum_{\substack{ J\subseteq J_{n+1}\\|J|=2s\\x_{n+1}\notin J}}
   \mkern -15mu P_{n+1}(J)Q(J)
   +
   \sum_{s=0}^{\left\lfloor\mkern-5mu\copy0\right\rfloor}\mkern-4mu(-1)^s2^{(n+1)-1-s}\mkern-15mu\sum_{\substack{ J\subseteq J_{n+1}\\|J|=2s\\x_{n+1}\in J}}
   \mkern -8mu P_{n+1}(J)Q(J),
 \]
 as required.
\end{proof}

 Note that as an immediate consequence of Lemma \ref{lemma: deal with
 Venkova's guy}, when $m$ is odd,
 $x_1\circ x_2\circ \cdots \circ x_{\frac{m+1}{2}}\in
 CP(\finitenonunitgrass{m})$ is in fact a central polynomial (that is,
 in $CP(\finitenonunitgrass{m})-T(\finitenonunitgrass{m})$). For if we
 evaluate by setting $x_i=e_{2i-1}e_{2i}$ for $i=1,2,\ldots,(m-1)/2$,
 and $x_{(m+1)/2}=e_{m}$, the result is $2^{(m-1)/2}e_1e_2\cdots e_m\ne0$.

\begin{lemma}\label{lemma: final touch}
 Let $m$ be an odd positive integer. Then each extremal element of 
 $\bss{m}$ is congruent modulo $U(m)$ to a linear combination of
 non-extremal elements of $\bss{m}$.
\end{lemma}

\setbox0=\hbox{$\scriptstyle\frac{m+1}{4}$}%\showthe\dp0
\dp0=2.5pt
\ht0=0pt
\wd0=15pt

\begin{proof}
 By Lemma \ref{lemma: deal with Venkova's guy}, 
 \[
  x_1\circ x_2\circ\cdots\circ x_{(m+1)/2}\cong\mkern-7mu
  \sum_{s=0}^{\left\lfloor\mkern-5mu\copy0\right\rfloor}\mkern-4mu(-1)^s2^{(m-1)/2-s}
  \mkern-25mu\sum_{\substack{ J\subseteq J_{(m+1)/2}\\|J|=2s}}
  \mkern -25mu P_{(m+1)/2}(J)Q(J)\ \mod{T^{(3)}},
 \]
 and since $x_1\circ x_2\circ\cdots\circ x_{(m+1)/2}\in S(m)$, it follows that
 \[
  P_{(m+1)/2}(\nullset)\cong\sum_{s=1}^{\left\lfloor\mkern-5mu\copy0\right\rfloor}
  (-1)^{s+1}2^{-s}\mkern-10mu\sum_{\substack{ J\subseteq J_{(m+1)/2}\\|J|=2s}}
  \mkern -8mu P_{(m+1)/2}(J)Q(J)\ \mod{U(m)}.
 \]
 Let $u\in \bss{m}$ be extremal, so $\deg(u)\le m$ and
 $\deg(u)-\lend{u}=(m+1)/2$, so $0\le\lend{u}\le (m-1)/2$. We claim that
 $u$ is obtainable from $P_{(m+1)/2}(\nullset)= x_1x_2\cdots
 x_{(m+1)/2}$ by substitution. If $\lend{u}=0$, this is immediate, so we
 may suppose that $s=\lend{u}>0$. If $\lbeg{u}=0$, then $u\in S\subseteq
 U(m)$, and in such a case, there is nothing to show. Thus we may
 suppose that $t=\lbeg{u}>0$ as well. Let $z=(m+1)/2-\lend{u}$, and for
 convenience, let $y_i=x_{z+i}$, $i=1,2,\ldots,\lend{u}$. Suppose that 
 $u=\prod_{r=1}^t x_{i_r}\prod_{r=1}^s\bigl(\com x_{j_{2r-1}},{x_{j_{2r}}}
 x_{j_{2r-1}}^{\beta_{2r-1}}x_{j_{2r}}^{\beta_{2r}}\bigr)\cong
  \prod_{r=1}^t x_{i_r}\prod_{r=1}^s\bigl(
 x_{j_{2r-1}}^{\beta_{2r-1}}x_{j_{2r}}^{\beta_{2r}}\bigr)\prod_{r=1}^s
 \com x_{j_{2r-1}},{x_{j_{2r}}}\mod{T^{(3)}}$. Since $z=\deg(u)-2s$, which
 is the number of copies of variables
 $x_{i_1},\ldots,x_{i_t},x_{j_1},\ldots,x_{j_{2s}}$, we may obtain
 $\prod_{r=1}^t x_{i_r}\prod_{r=1}^s\bigl(
 x_{j_{2r-1}}^{\beta_{2r-1}}x_{j_{2r}}^{\beta_{2r}}\bigr)$ from
 $x_1x_2\cdots,x_z$ by substitution. Then let $y_i$ be replaced by $\com
 x_{2i-1},{x_{2i}}$ for $i=1,2,\ldots, s$, with $u$ as the end result.
 After  substitution in the summation that is congruent to
 $P_{(m+1)/2}(\nullset)$ modulo $U(m)$, any summand that has a variable
 $y_i$ appearing in $Q(J)$ will be $0\mod{T^{(3)}}$, so the only terms that
 survive the substitution are those that have every $y_i$ appearing in
 $P_{(m+1)/2}(J)$. Thus each surviving summand contains every commutator
 factor of $u$, plus those from $J$. Since $J\ne\nullset$ for each
 summand, this means that each summand has greater end length than $u$.
 However, since each summand was (before substitution) multilinear of
 degree $(m+1)/2$, it follows that each summand after summation has the
 same degree as $u$, and so no summand is extremal, as required. By
 application of the various identities in Lemma \ref{lemma: handy}, each
 summand can be converted, modulo $T^{(3)}$, into one or more elements of
 $\bss{m}$ of the same degree, each with equal or greater end length
 than that of the summand. The result follows now. 
\end{proof}

\begin{definition}\label{definition: def of M}
 Let $t$ and $n$ be positive integers with $t\le n$, and define $M_{t,n}$ to be the set of all elements of the form 
 \begin{gather}
  x_t^{\alpha_t} \prod_{r=1}^s \com x_{j_{2r-1}},{x_{j_{2r}}}x_{j_{2r-1}}^{\beta_{2r-1}}x_{j_{2r}}^{\beta_{2r}} \label{Mdef i}\\ 
 \intertext{or}
  \biggl( \prod_{r=1}^l x_{i_r}^{\alpha_r}\biggr)x_t^{\alpha_t}\prod_{r=1}^s \com x_{j_{2r-1}},{x_{j_{2r}}}x_{j_{2r-1}}^{\beta_{2r-1}}x_{j_{2r}}^{\beta_{2r}}, \label{Mdef ii} \\
 \intertext{or, in the case $t=n$ only, elements of the form}
 \prod_{r=1}^n x_{i_r}^{\alpha_r} \label{Mdef iii} 
 \end{gather}
 as well, subject to the following requirements.
 \begin{list}{(\alph{parts})}{\usecounter{parts}} 
 \item\label{M (a)} % (i)
  $s,l\ge1$.
 \item % (ii)  
  $1 \le i_1<\cdots<i_{l}<t$ for elements of type (3.2) or (3.3).
  \item % (iii)
   $j_1<\cdots<j_{2s}$.
  \item % (iv) 
   $\{i_1,\ldots,i_l,t\} \cap \{j_1,\ldots,j_{2s}\}=\nullset$.
  \item % (v) 
   $\{i_1,\ldots,i_l,t\} \cup\{j_1,\ldots,j_{2s}\}=\{1,\ldots,n\}$.
  \item % (vi) 
    For every $1 \le k \le l$, $1\le \alpha_m,\alpha_k\le \frac{m+1}{2}$, and, if $p>2$, $\alpha_m,\alpha_k\le p-1$.
  \item % (vii)    
     For $1 \le r \le 2s$, $0\le\beta_r\le \frac{m+1}{2}$, and, if $p>2$, $\beta_r \le p-1$.
 %% fixed since first submission, needed the upper bound of (m+1)/2
 \end{list}
\end{definition}

\begin{lemma}\label{lemma: m fact} 
  Let $t,n$ be positive integers with $t\le n$. Then for any $u \in M_{t,n}$, there exist
  $g_1,g_2,\ldots,g_n$ in $\finitenonunitgrass{z}$, where $z=2(\deg(u)-\lend{u})-1$, such that the following hold:
   \begin{list}{(\roman{parts})}{\usecounter{parts}}
   \item $0\ne u(g_1,g_2,\ldots,g_n)\in \gen{\prod_{r=1}^z e_i}$;
   \item for any $v\in M_{t,n}$ with $u>v$, $v(g_1,g_2,\ldots,g_n)=0$;
   \item if $p>2$, then for any $g\in \nonunitgrass$, $\com g_t,g g_t^{p-1}=0$.
  \end{list}
\end{lemma}

\begin{proof}
%% fixed since submission, was Lemma 2.6
 This result follows immediately from Lemma 2.7 of \cite{Ra}, upon noting that $u\in M_{t,n}$ implies that
 every factor of the form $x^r$ for $x\in X$ has $r< p$, so $r!\not\cong 0\mod{p}$, and $r<\frac{m+1}{2}$,
 so the evaluation of $x^r$ in the lemma is always possible.
\end{proof}

\begin{theorem}\label{theorem: central polys for finite nonunitary grassman}
 For any positive integer $m$, $CP(\finitenonunitgrass{m})=U(m)$.
\end{theorem}

\begin{proof}
 By \cite{Ra}, Theorem \thmMainTheoremForNonunitary, 
 $S\subseteq CP(\nonunitgrass)$, and $CP(\nonunitgrass)\subseteq 
 CP(\finitenonunitgrass{m})$, so $S\subseteq CP(\finitenonunitgrass{m})$. 
 If $m$ is even, $S(m)=S$ and thus $S(m)\subseteq CP(\finitenonunitgrass{m})$, while if 
 $m$ is odd, then by Corollary \ref{corollary: new finite nonunitary central guy},
 $S(m)=S+\set x_1\circ\cdots \circ x_{\frac{m+1}{2}}\endset^S \subseteq CP(\finitenonunitgrass{m})$.
 Thus $U(m)=S(m)+T(\finitenonunitgrass{m})\subseteq CP(\finitenonunitgrass{m})$.

 For the converse, we note that $CP(\finitenonunitgrass{m})$ is a
 $T$-space, and thus by Lemma \ref{lemma: basic essential} (ii),
 $CP(\finitenonunitgrass{m})$ is the linear span of its essential
 polynomials. It suffices therefore to prove that any essential element
 of $CP(\finitenonunitgrass{m})$ belongs to $U(m)$. So let $f\in
 CP(\finitenonunitgrass{m})$  be essential in
 $CP(\finitenonunitgrass{m})$, and suppose that $f\notin U(m)$. By Lemma
 \ref{lemma: at, lemma 6} if $m$ is even, or by Lemma \ref{lemma: at,
 lemma 11} if $m$ is odd, $f$ is congruent modulo $T(G_m)$ to a linear
 combination $\sum a_i u_i$ of $u_i\in \bss{m}$. By Lemma \ref{lemma:
 final touch}, we may assume that no $u_i$ is extremal, so in the latter
 case, we would have $2\deg(u_i)-2\lend{u_i}\le m-1$. Without loss of
 generality, we may assume that the variables that appear in $f$ are
 $x_1,x_2,\ldots,x_n$ for some positive integer $n$. We may further
 assume that for each $i$, $u_i$ is essential in $\kzerox$ in the
 variables $x_1,x_2,\ldots,x_n$. For suppose to the contrary that for
 some $i$ and $j$, $x_j$ does not appear in $u_i$. Let
 $f_0=\sum_{x_j\text{ not in }u_r} a_ru_r$ and $f_0=\sum_{x_j\text{in
 }u_r} a_ru_r$, so $f\cong f_0+f_1\mod{T(\finitenonunitgrass{m})}$.
 Since $x_j$ appears in every monomial in $f$ and in $f_1$, we have
 $0\cong f\rest{x_j=0}=f_0+(f_1\rest{x_j=0})= f_0\mod{U(m)}$, and so
 $f\cong f_1\mod{U(m)}$.

 Next, we observe that for any $i$, if $\lbeg{u_i}=0$, then by \cite{Ra}, Corollary \corAltDescriptionOfSGen, $u_i\in S+T^{(3)}\subseteq 
 T(\finitenonunitgrass{m})$. Thus for each $i$, $\lbeg{u_i}>0$. Let  
 \begin{align*}
   d=\max\{\, j\mid 1\le j\le n\ \text{ and }&\text{there exists $i$ such that}\\ 
  &\text{$x_j$ appears in the beginning of $u_i$}\,\}.
 \end{align*}
 Let $X=\set i\mid x_d \text{ appears in the beginning of
 $u_i$}\endset$, and note that for each $i\in X$, $u_i\in M_{d,n}$.  Let
 $f_d=\sum_{i\in X} a_i u_i$, so that $f_d\in \gen{M_{d,n}}$, and let
 $f_e=\sum_{i\notin X} a_iu_i$, so that $f\cong f_d+f_e\mod{U(m)}$.

 Suppose that $f_e\in U(m)$. Then $f_d\cong f\not\cong 0\mod{U(m)}$, and
 so $f_d\in \gen{M_{d,n}}-\set 0\endset$. Let $j\in X$ be such that
 $u_j>u_i$ for all $i\in X$, $i\ne j$. By Lemma \ref{lemma: m fact},
 there exist $g_1,g_2,\ldots,g_n\in \finitenonunitgrass{m}$ such that
 $f_d(g_1,g_2,\ldots,g_n)$ is a scalar multiple of $e_1e_2\cdots e_z$,
 where $z=2\deg(u_j)-2\lend{u_j}-1\le m-1$. Thus
 $f_d(g_1,g_2,\ldots,g_n)\notin C_{\finitenonunitgrass{m}}$ and so
 $f_d\notin CP(\finitenonunitgrass{m})$, which means that $f\notin
 CP(\finitenonunitgrass{m})$. Since this is not the case, it follows
 that $f_e\notin U(m)$.  Consider $j$ such that  $\alpha_j u_j$ is a
 summand of $f_e$. Suppose that $\deg_{x_d} u_j<p$, so there exist
 $0<i_1<i_2<\cdots <i_t<d$, positive integers
 $\alpha_1,\ldots,\alpha_t$, $0<j_{1}<\cdots<j_{2s}$, nonnegative
 integers $\beta_1,\ldots,\beta_{2s}$ such that $m\in\set
 j_1,\ldots,j_{2s}\endset$, $\beta_d\le p-2$, and
 \[
  u_j=\prod_{k=1}^t x_{i_k}^{\alpha_k}\,\prod_{r=1}^s\com x_{j_{2r-1}},{x_{j_{2r}}} x_{j_{2r-1}}^{\beta_{2r-1}}x_{j_{2r}}^{\beta_{2r}}.
 \]

 We may assume without loss of generality that $d$ is odd, since the
 argument for the case when $d$ is even can be converted,  modulo $T^{(3)}$,
 to the case where $d$ is odd by a sign change in $\alpha_j$. Thus
 $d=2b-1$ for some $b$ with $1\le b\le s$, and for convenience, let
 $\beta=\beta_{j_{2b-1}}$. We apply \cite{Ra}, Corollary \corNeedForRGen, taking $v=\prod_{k=1}^t x_{i_k}^{\alpha_k}$ and 
 $u=\prod_{\substack{ r=1\\2r-1\ne d}}^s \com x_{j_{2r-1}},{x_{j_{2r}}} 
 x_{j_{2r-1}}^{\beta_{2r-1}}x_{j_{2r}}^{\beta_{2r}}$, to find that
 \begin{align*}
  u_j&\overset{S+T^{(3)}}{\cong} (\beta+1)^{-1}\com x_{j_{2d}},v x_d^{\beta+1}x_{j_{2b}}^{\beta_{2b}}u\\
    &\overset{T^{(3)}}{\cong}(\beta+1)^{-1}\sum_{k=1}^t\alpha_k\com x_{j_{2d}},{x_{i_k}} x_{i_k}^{\alpha_k-1}\bigl(\prod_{\substack{l=1\\l\ne k}}^tx_{i_l}^{\alpha_l}\bigr)x_d^{\beta+1}x_{j_{2b}}^{\beta_{2b}}u\quad%\text{by Lemma \ref{lemma: handy} (iii) and (vi)}\\
    \text{\begin{tabular}[t]{l} by Lemma \ref{lemma: handy}\\ (iii) and (vi)\end{tabular}}\\
    &\overset{T^{(3)}}{\cong}\sum_{k=1}^t(\beta+1)^{-1}\alpha_k\bigl(\prod_{\substack{l=1\\l\ne k}}^tx_{i_l}^{\alpha_l}\bigr)x_d^{\beta+1}\com x_{j_{2b}},{x_{i_k}} x_{j_{2b}}^{\beta_{2b}}x_{i_k}^{\alpha_k-1}u\quad\text{by Lemma \ref{lemma: handy} (vi)}\\
 \end{align*}
 Now for each $k$, $\displaystyle\bigl(\prod_{\substack{ l=1\\l\ne k}}^tx_{i_l}^{\alpha_l}\bigr)x_d^{\beta+1}
 \com x_{j_{2b}},{x_{i_k}} x_{j_{2b}}^{\beta_{2b}}x_{i_k}^{\alpha_k-1}u$
 is equal to
 \[
   \bigl(\prod_{\substack{ l=1\\l\ne k}}^tx_{i_l}^{\alpha_l}\bigr)x_d^{\beta+1}\com x_{j_{2b}},{x_{i_k}} x_{j_{2b}}^{\beta_{2b}}x_{i_k}^{\alpha_k-1}
   \prod_{\substack{ r=1\\r\ne b}}^s
             \com x_{j_{2r-1}},{x_{j_{2r}}} x_{j_{2r-1}}^{\beta_{2r-1}}x_{j_{2r}}^{\beta_{2r}},
 \]
 and by Lemma \ref{lemma: handy}, working modulo $T^{(3)}$, the end of this element
 can be rearranged so as to give an element of $R$ with beginning 
 $\bigl(\prod_{\substack{l=1\\l\ne k}}^tx_{i_l}^{\alpha_l}\bigr)x_d^{\beta+1}$.	     
 We have proven now that if $\deg_{x_d} u_j<p$, then $u_j$ is congruent to a linear 
 combination of elements of $M_{d,n}$. It follows that $d$ is such that there exist 
 $f_d\in \gen{M_{d,n}}$ and $f_e$ in
 \[
 \gen{\set u\in R\mid \deg_{x_d}u=p\ \text{and for every $i$, $i<d$ if $x_i$ 
  is in the beginning of $u$}\endset}
 \]
 such that $f\cong f_d+f_e\mod{(S+T^{(3)}}$. Suppose now that $d$ is minimal
 with respect to this property. If $f_d\cong 0\mod{S+T^{(3)}}$,  then
 $f\cong f_e\mod{S+T^{(3)}}$, which contradicts our choice of $d$ since
 $f_e$ is a linear combination of elements of $R$ in whose beginning 
 only elements $x_i$ with $i<d$ appear. Thus $f_d\in
 \gen{M_{d,n}}-\set0\endset$, and so by Lemma \ref{lemma: m fact}, there
 exist $g_1,g_2,\ldots,g_n\in \finitenonunitgrass{m}$ such that (as
 above) $f_d(g_1,g_2,\ldots,g_n)\notin C_{\finitenonunitgrass{m}}$ and
 for any $g\in \finitenonunitgrass{m}$ and any positive integer $\beta$,
 $\com g_m,g g_m^{p-1}g^\beta=0$. Since $f_e=\sum \gamma_iv_i$, where
 for each $i$, $v_i\in R$ and $\deg_{x_d}v_i=p$, it follows that in each
 $v_i$, $x_d$ appears in a term of the form $\com
 x_d,gx_d^{p-1}g^{\beta}$, and so $v_i(g_1,g_2,\ldots,g_n)=0$ for each
 $i$. Thus $f(g_1,g_2,\ldots,g_n)\cong
 f_d(g_1,g_2,\ldots,g_n)+0\mod{C_{\finitenonunitgrass{m}}}$. But
 $f_d(g_1,g_2,\ldots,g_n)\notin C_{\finitenonunitgrass{m}}$, and thus
 $f\notin CP(\finitenonunitgrass{m})$. Since this contradiction follows
 from our assumption that $f\notin U(m)$, it follows that $f\in U(m)$,
 as required.
\end{proof}

We note that by \cite{At} Lemma 9, if $p>2$, then $w_n\in T(\finitenonunitgrass{m})$ for
$n\ge \frac{m+1}{2(2p-1)}$. We shall let $r_0=\frac{m+1}{2(2p-1)}$.

\begin{corollary}\label{corollary: t space basis for finite nonunitary}
 If $k$ if a field of characteristic zero, then
 \[
   CP(\finitenonunitgrass{m})=\set \com x_1,{x_2}, \com x_1,{x_2}\com x_3,{x_4},x_1\circ\cdots \circ x_{b_m},
   (x_1\circ\cdots \circ x_{b_m})x_{b_m+1}\endset^S
   \] 
  where $b_m=\left\lfloor\frac{m}{2}\right\rfloor$, while if $k$ is a field of characteristic $p>2$, then 
 \begin{align*}
   CP(\finitenonunitgrass{m})&=\{\, \com x_1,{x_2}, \com x_1,{x_2}\com x_3,{x_4},x_1\circ\cdots \circ x_{b_m},
  (x_1\circ\cdots \circ x_{b_m})x_{b_m+1}\,\}^S\\
  &\hskip-10pt+\{\,x_1^p,\,x_2x_1^p\,\}^S+\{\, w_k\mid 1\le k\le r_0\,\}^S\\
  &+\{\, u\mid \text{if $m$ is odd and $r_0\in\mathbb{Z}$, then $u=x_{2r_0+1}w_{r_0}$, else $u=0$}\,\}^S
 \end{align*}
\end{corollary}

\end{document}